\documentclass[12pt]{amsart}
\usepackage{amssymb}
\usepackage{amsmath,amssymb,amsfonts,amsthm,graphics,
latexsym, amscd, amsfonts, epsfig,
eepic,epic,psfrag,setspace,eufrak}
\usepackage{e-jc}
\usepackage{mathrsfs}
\usepackage{color}
\usepackage{eucal}
\usepackage{eufrak}
\usepackage[all]{xypic}
\usepackage{xspace}
\usepackage{algorithm}
\usepackage{algorithmic}

\theoremstyle{plain}
\newtheorem{theorem}{Theorem}
\newtheorem{corollary}{Corollary}

\newtheorem{lemma}{Lemma}

\theoremstyle{definition}

\theoremstyle{remark}

\numberwithin{equation}{section}

\begin{document}

\begin{center}
{\bf\Large On the uniform generation of modular diagrams} \\
\vspace{15pt} Fenix W.D. Huang and Christian M. Reidys$^{\,\star}$
\end{center}

\begin{center}
Center for Combinatorics, LPMC-TJKLC\\
Nankai University  \\
Tianjin 300071\\
         P.R.~China\\
         Phone: *86-22-2350-6800\\
         Fax:   *86-22-2350-9272\\
$^{\,\star}$duck@santafe.edu
\end{center}

\centerline{\bf Abstract}
In this paper we present an algorithm that generates $k$-noncrossing,
$\sigma$-modular diagrams with uniform probability.
A diagram is a labeled graph of degree $\le 1$ over $n$ vertices drawn
in a horizontal line with arcs $(i,j)$ in the upper half-plane.
A $k$-crossing in a diagram is a set of
$k$ distinct arcs $(i_1, j_1), (i_2, j_2),\ldots,(i_k, j_k)$ with
the property $i_1 < i_2 < \ldots < i_k  < j_1 < j_2 < \ldots< j_k$.
A diagram without any $k$-crossings is called a $k$-noncrossing diagram
and a stack of length $\sigma$ is a maximal sequence
$((i,j),(i+1,j-1),\dots,(i+(\sigma-1),j-(\sigma-1)))$.
A diagram is $\sigma$-modular if any arc is contained in a stack of length
at least $\sigma$.
Our algorithm generates after $O(n^k)$ preprocessing time,
 $k$-noncrossing, $\sigma$-modular diagrams in $O(n)$ time
and space complexity.


\section{Introduction}\label{S:Intro}

A ribonucleic acid (RNA) molecule is the helical configuration of a primary
structure of nucleotides, {\bf A}, {\bf G}, {\bf U} and {\bf C}, together
with Watson-Crick ({\bf A}-{\bf U}, {\bf G}-{\bf C}) and ({\bf U}-{\bf
G}) base pairs (arcs).
It is well-known that RNA structures exhibit cross-serial nucleotide
interactions, called pseudoknots. First recognized in the turnip yellow
mosaic virus in \cite{Rietveld:82}, they are now known to be widely
conserved in functional RNA molecules.

Modular $k$-noncrossing diagrams represent a model of RNA pseudoknot
structures \cite{Reidys:07pseu,Reidys:07lego}, that is RNA
structures exhibiting cross-serial base pairings. The particular case
of modular noncrossing diagrams, i.e.~RNA secondary structures have been
extensively studied \cite{Schuster:98,Waterman:93,Waterman:78b,Waterman:79}.

A diagram is a labeled graph over the vertex set $[n]=\{1, \dots, n\}$ with
vertex degrees not greater than one. The standard representation of a
diagram is derived by drawing its vertices in a horizontal line and
its arcs $(i,j)$ in the upper half-plane. A $k$-crossing is a set of
$k$ distinct arcs $(i_1, j_1), (i_2, j_2),\ldots,(i_k, j_k)$ with
the property
\begin{equation*}
i_1 < i_2 < \ldots < i_k  < j_1 < j_2 < \ldots< j_k.
\end{equation*}
A diagram without any $k$-crossings is called a $k$-noncrossing diagram.
Furthermore, a stack of length $\sigma$ is a maximal sequence
of ``parallel'' arcs,
\begin{equation*}
((i,j),(i+1,j-1),\dots,(i+(\sigma-1),j-(\sigma-1)))
\end{equation*}
and is also referred to as a $\sigma$-stack.
A $k$-noncrossing diagram having only stacks of lengths one is called a core.
\begin{figure}[ht]
\centerline{\epsfig{file=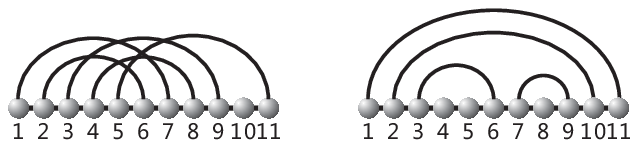,width=0.8\textwidth} \hskip8pt}
\caption{\small $k$-noncrossing diagrams: a $4$-noncrossing diagram (left) and
a $2$-noncrossing diagram (right). The arcs $(2,6), (4,8)$ and $(5,11)$ form
a $3$-crossing in the left diagram.
} \label{F:k-non}
\end{figure}

Biophysical structures do not exhibit any isolated bonds. That is, any
arc in their diagram representation is contained in a stack of length at
least two. We call a diagram, whose arcs are contained in stacks of lengths
at least $\sigma$, $\sigma$-modular. Modular, $k$-noncrossing
diagrams are likely candidates for natural molecular structures.
Sequence lengths of interest for such structures range from $75$--$300$
nucleotides.

\begin{figure}[ht]
\centerline{\epsfig{file=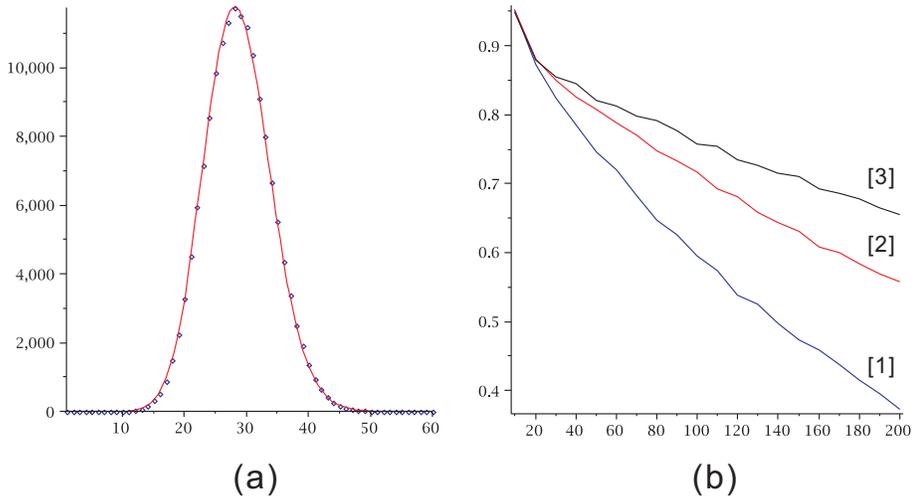,width=0.8\textwidth} \hskip8pt}
\caption{\small Uniformity and success-rate of Algorithm~$2$.
We run Algorithm~$2$ for $5\times 10^6$ times attempting to generate
$3$-noncrossing $2$-modular diagrams over $20$ vertices. $4,354,410$ of these
executions generate a modular diagram.
In (a) we display the frequency distribution of multiplicities (dots)
and the Binomial distribution (curve).
In (b) we display the success rate of Algorithm~$2$ as a function
of $n$ for the following classes of modular diagrams: $k=3$, $\sigma=2$
([1]), $k=4$, $\sigma=2$ ([2]) and $k=5$, $\sigma=2$ ([3]).
} \label{F:dis_suc}
\end{figure}

The main result of this paper is an algorithm that generates $k$-noncrossing,
$\sigma$-modular diagrams with uniform
probability. Our construction is motivated by the ideas of \cite{PNAS:09}, where
a combinatorial algorithm has been presented that uniformly generates
$k$-noncrossing diagrams in $O(n^k)$ time complexity. To be precise, we
generate $k$-noncrossing modular diagrams ``locally'' having a
success rate that depends on specific parameters, see Fig.~\ref{F:dis_suc}.

The paper is organized in two sections. In Section~\ref{S:core} we lay the
foundations for our main result by generating core diagrams with uniform
probability. In Section~\ref{S:can} we introduce weighted cores and subsequently
prove the main theorem.


\section{Core diagrams}\label{S:core}

\begin{figure}[ht]
\centerline{\epsfig{file=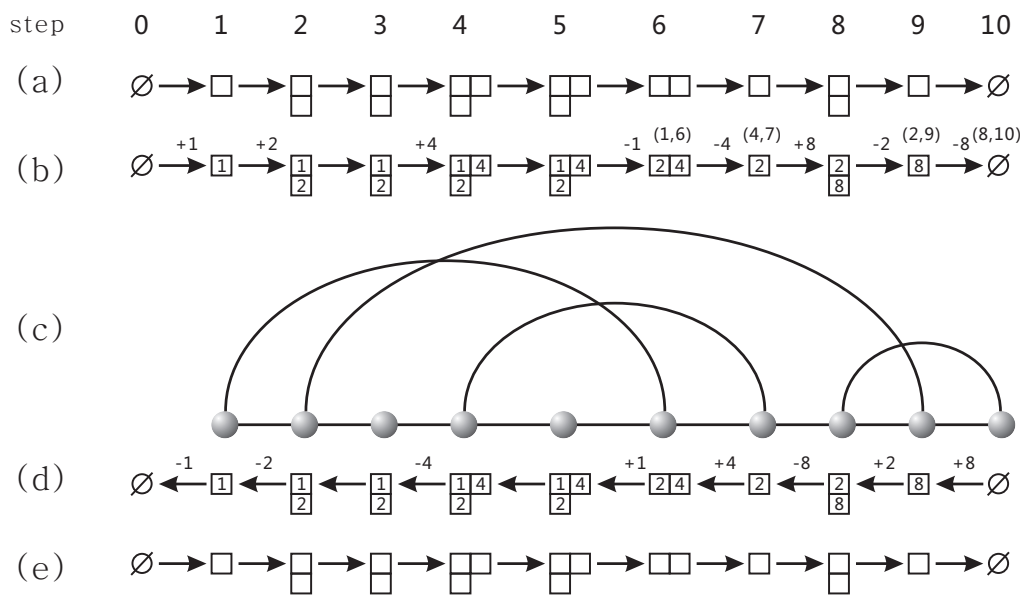,width=0.9\textwidth} \hskip8pt}
\caption{\small From $*$-tableaux to diagrams and back. 
Reading (a) from left to right, we insert $i$ into the new square in case of
$\lambda^i \setminus \lambda^{i-1}$ being a $-\square$-step and extract 
the square via inverse RSK if $\lambda^i \setminus \lambda^{i-1}$ is a 
$-\square$-step. The extraction leads to an arc. 
Reading (c) from right to left, $\lambda^{i-1}$ is
obtained by RSK insertion of $j$ into $\lambda^i$ if $i$ is the terminal
of an arc. We do nothing if $i$ is an isolated vertex and we remove the
square with entry $i$ in case of $i$ being an origin of an arc.} \label{F:bijection}
\end{figure}

A shape $\lambda$ is a set of squares
arranged in left-justified rows with weakly decreasing number of boxes in each
row. A Young tableau is a filling in squares in the shape with numbers, which is
weakly increasing in each row and strictly increasing in each column. A
$*$-tableau of $\lambda^n$ is a sequence of shapes,
$$
\varnothing=\lambda^0, \lambda^1, \ldots, \lambda^n,
$$
such that $\lambda^i$ is differ from $\lambda^{i-1}$ by at most one
square. See Fig.~\ref{F:bijection} (a).

According to \cite{Chen} we have a bijection between $k$-noncrossing diagrams
and a $*$-tableaux of $\varnothing$ having at most $(k-1)$ rows.
Let us make the bijection explicit: reading the $*$-tableaux having $n$
steps from left to right we do the following:
if $\lambda^i\setminus\lambda^{i-1}=+\square$, we insert $i$ in the new
square. Otherwise if $\lambda^i\setminus \lambda^{i-1}=-\square$, we
extract the unique entry $j$ via inverse RSK algorithm \cite{Chen,PNAS:09} and
form an arc $(j,i)$. By inverse RSK algorithm we mean the following:
given a Young tableau $Y^i$ of shape $\lambda^i$ and a shape $\lambda^{i+1}$ such
that $\lambda^{i+1}\setminus \lambda^i=-\square$, there exists a unique entry $j$ 
of $Y^i$ and a Young tableau $Y^{i+1}$ of shape $\lambda^{i+1}$ such that 
RSK-insertion of $j$ into $Y^{i+1}$ recovers $Y^i$.
Finally, in case of $\lambda^i\setminus \lambda^{i-1}=\varnothing$ we do nothing,
see Fig.~\ref{F:bijection}. Given a $k$-noncrossing diagram, we read the
vertices from right to left and initialize $\lambda^n=\varnothing$. If $i$ is a
terminal of an arc, $(j,i)$, we obtain $\lambda^{i-1}$ by inserting $j$ into
$\lambda^i$ via RSK insertion. If $i$ is an isolated vertex we do nothing, and
remove the square contain $i$ when it is an origin of an arc, see
Fig.~\ref{F:bijection}.

Let
\begin{equation*}
T_i(\lambda)=\{ (\lambda^j)_{0\le j\le i} \mid \text{\rm
$(\lambda^j)_{j}$ is a $*$-tableau having at most $(k-1)$ rows and
$\lambda^i=\lambda$}\}.
\end{equation*}
Any $\vartheta \in T_i(\lambda)$ induces a unique
arc-set $A(\vartheta)$.
We set $A_0(\vartheta)=\varnothing$ and do the following
\begin{itemize}
 \item for a $+\square$-step, we insert $h$ into the new square,
 \item for a $\varnothing$-step, we do nothing,
 \item for a $-\square$-step, we extract the unique entry, $j(h)$, of the
tableaux $Y^{h-1}$ which, if RSK-inserted into $Y^h$, recovers $Y^{h-1}$
and set $A_h(\vartheta)=A_{h-1}(\vartheta)\dot\cup\{(j(h),h)\}$.
\end{itemize}
Setting $A(\vartheta)=A_i(\vartheta)$ we obtain an induced arc set
$A(\vartheta)$, as well as a unique sequence of Young tableaux
$Y(\vartheta)=\{Y^0=\varnothing,Y^1,\ldots, Y^i\}$, where for $h\le i$, $Y^h$
is a Young tableau of shape $\lambda^h$. These extractions generate a set of
arcs $(j(i),i)$, which in turn uniquely determines a $k$-noncrossing
diagram.

\begin{lemma} \label{L:removal}
Suppose $r\ge 1$ and $\vartheta_{p,q,r}\in T_i(\lambda)$ is a $*$-tableaux
such that
\begin{equation*}
(p,q),(p+1,q-1),\dots,(p+r,q-r)
\end{equation*}
are stacked pairs of insertion-extraction steps.
Let $f(\vartheta_{p,q,r})\in T_i(\lambda)$ be the $*$-tableaux in which
all $r$ insertion-extraction pairs $(p+1,q-1),\dots,(p+r,q-r)$ are replaced by
$2r$ $\varnothing$-steps.
Then we have a correspondence between $\vartheta_{p,q,r}$ and
$f(\vartheta_{p,q,r})$.
\end{lemma}
\proof Let $Y(\vartheta_{p,q,r})$ denote its associated sequence of Young tableaux,
\begin{equation}
(Y^t)_{0\le t\le i}=(Y^0=\varnothing, Y^1, \ldots, Y^i).
\end{equation}
We next construct a new sequence of Young tableaux,
\begin{equation}
Y(f(\vartheta_{p,q}))=\{J^0,J^1,\ldots, J^n=Y^i\},
\end{equation}
from right to left via the following algorithm
\begin{itemize}
\item for a $-\square$-step of the original $*$-tableaux, $\vartheta_{p,q,r}$,
      let $j$ be the unique entry extracted from $Y^{t-1}$ which if
      RSK-inserted into $Y^t$ recovers $Y^{t-1}$.
      If $t=q,q-1,\dots,q-r$ we do nothing, otherwise: $J^{t-1}$ is obtained by
      RSK-insertion of $j$ into $J^t$,
\item for a $\varnothing$-step, we do nothing,
\item for a $+\square$-step, if $t=p+1,\dots,p+r$, we do nothing, otherwise
      $J^{t-1}$ is obtained by removing the square with entry $t$ from $J^t$.
\end{itemize}
By construction, $J^0=\varnothing$ and considering the induced sequence of
shapes of the sequence of Young tableaux $J^0,\ldots, J^i$ we obtain a
unique $*$-tableau $f(\vartheta_{p,q,r})$.
By construction $f(\vartheta_{p,q,r})$ has $\varnothing$-steps at step
$p+1,\dots,p+r$ and steps $q-1,\dots,q-r$, respectively.

\begin{figure}[ht]
\centerline{\epsfig{file=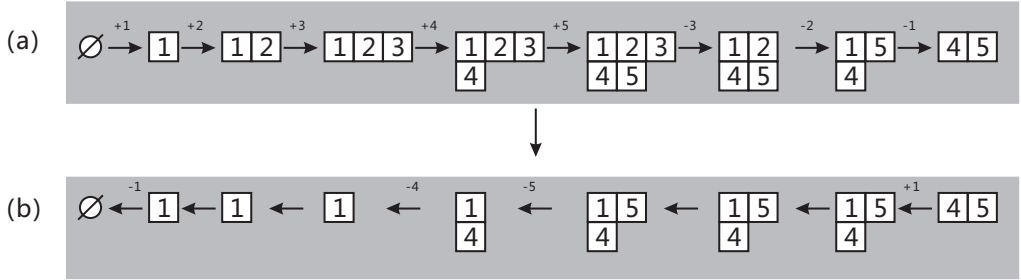,width=0.9\textwidth} \hskip8pt}
\caption{\small (a) a $*$-tableaux $\vartheta_{1,8,2}$ in which $(1,8),
(2,7)$ and $(3,6)$ are stacked pairs of insertion-extraction
steps. (b) $f(\vartheta_{1,8,2})$ is the unique $*$-tableau derived from
$\vartheta_{1,8,2}$ in which the $4$ steps: step $2,3,6$ and  $7$ are
$\varnothing$-steps. } \label{F:path}
\end{figure}

Suppose we are given a $*$-tableaux $\psi_{p,q,r}$ having the
insertion-extraction pair $(p,q)$ and $\varnothing$-steps at step $p+1,
\dots,p+r$ and $q-1,\dots,q-r$, respectively together with its sequence
of Young tableaux $(J^t)_{0\le t \le i}$.  Then we construct the sequence
of Young tableaux $(Y^t)_{0\le t \le i}$ initialized $Y^0=J^0=\varnothing$:
\begin{itemize}
\item for a $-\square$-step of the original $*$-tableaux, $\psi_{p,q,r}$,
      let $j$ be the unique entry extracted from $Y^{t-1}$ which if
      RSK-inserted into $Y^t$ recovers $Y^{t-1}$.
      $Y^{t-1}$ is obtained by RSK-insertion of $j$ into $Y^t$,
\item for a $\varnothing$-step of $\psi_{p,q,r}$, if $t=q-1,\dots,q-r$,
      we add a square and insert $p+1,\dots,p+r$.
      If $t=p+1,\dots,p+r$, we remove the square with
      the respective entry $p+1,\dots,p+r$.
      Otherwise, we do nothing.
\item for a $+\square$-step of $\psi_{p,q,r}$, $Y^{t-1}$ is obtained by removing
      the square with entry $t$.
\end{itemize}
It is straightforward to verify that the above algorithm is welldefined and
recovers the $*$-tableaux $\vartheta_{p,q,r}$ from
$f(\vartheta_{p,q,r})$, whence the lemma. See Fig.~\ref{F:path}.
\qed

We next consider
\begin{equation}
T^c_i(\lambda)= \{t \in T_i(\lambda)\mid \forall a\in A(t), a
\text{ is an isolated arc} \}
\end{equation}
and set $t^c_i(\lambda)=|T^c_i(\lambda)|$.
Given a shape $\lambda^i$, let $\lambda_{j+}^{i-1}$ denote the shape
from which $\lambda^i$ is obtained by adding a square in the $j$th row, and
$\lambda_{j-}^{i-1}$ denote the shape from which $\lambda^i$ is derived by
removing a square in the $j$th row. Thus tracing back a shape $\lambda^i$
we observe that it is either derived by
\begin{itemize}
\item $\lambda^{i-1}_{j+}$ (obtained by adding a square in the $j$th row),
\item $\lambda_0^{i-1}$ (doing nothing), or
\item $\lambda^{i-1}_{j-}$ (obtained by removing a square in the $j$th row).
\end{itemize}
\begin{lemma}\label{L:recursion}
\begin{equation}  \label{E:rec_c}
 t^c_i(\lambda^i)=t^c_{i-1}(\lambda_0^{i-1})+\sum_{j=1}^{k-1}
t^c_{i-1}(\lambda^{i-1}_{j+})+
\sum_{j=1}^{k-1} \sum_{p=0}^{\lfloor \frac{i-1}{2}
\rfloor}(-1)^p t^c_{i-1-2p}(\lambda^{i-1-2p}_{j-}).
\end{equation}
\end{lemma}
\proof
By construction, $+\square$-steps as well as $\varnothing$-steps do not
induce new arcs. An arc $\alpha$ is only formed when removing a square
and such an arc is potentially stacking. Let
\begin{equation*}
G_{i-1}(\lambda_{j-}^{i-1})=\{(\lambda^h)_{0\le h\le i-1} \in
T^c_{i-1}(\lambda_{j-}^{i-1})\mid
\text{\rm $\lambda^{i}\setminus\lambda^{i-1}=
-\square_j$ and $\alpha$ is stacking} \}.
\end{equation*}
Thus, for any $t\in T_{i-1}^c(\lambda_{j-}^{i-1})\setminus G_{i-1}(
\lambda_{j-}^{i-1})$, the $*$-tableaux $(t,-\square_j)$ is contained in
$T_i^c(\lambda^i)$. We accordingly arrive at
\begin{equation}
T_i^c(\lambda)=T_{i-1}^c(\lambda_0^{i-1}) \dot\cup
 \left(\bigcup_{j=1}^{k-1} T_{i-1}^c(\lambda^{i-1}_{j+})\right)
 \dot\cup \left(\bigcup_{j=1}^{k-1}[T_{i-1}^c(\lambda^{i-1}_{j-})\setminus
                G_{i-1}(\lambda^{i-1}_{j-})]\right)
\end{equation}
which implies
\begin{equation} \label {E:union}
t^c_i(\lambda^i)= t^c_{i-1}(\lambda^{i-1}_0)+
             \sum_{j=1}^{k-1}t^c_{i-1}(\lambda^{i-1}_{j+})+
             \sum_{j=1}^{k-1}\left[t^c_{i-1}(\lambda^{i-1}_{j-})-
             g_{i-1}(\lambda^{i-1}_{j-})\right].
\end{equation}
We next provide an interpretation of $G_{i-1}(\lambda^{i-1}_{j-})$.
Suppose the entry extracted at step $i$ is $j(i)$. The fact that $\alpha$
is in a stack implies that the $(i-1)$th step is also a $-\square$ step
and that the extracted entry is $j(i)+1$. For $\vartheta\in
G_{i-1}(\lambda^{i-1}_{j-})$, we apply Lemma~\ref{L:removal} and replace the
insertion of step $j(i)+1$ and the extraction at step $(i-1)$ by
respective $\varnothing$-steps, and thereby obtain the $*$-tableaux
$f(\vartheta)$. We then remove the two $\varnothing$-steps and
obtain the unique $*$-tableaux
\begin{equation*}
\vartheta^\prime \in T^c_{i-3}(\lambda^{i-3}_{j-}),
\end{equation*}
where $\lambda^i$ can be derived from $\lambda^{i-3}_{j-}$ by removing a
square in the $j$th row.
We next claim $\vartheta^\prime \in T^c_{i-3}(\lambda^{i-3}_{j-})\setminus
G_{i-3}(\lambda^{i-3}_{j-})$. Suppose $\vartheta^\prime \in G_{i-3}(
\lambda^{i-3}_{j-})$, then $\vartheta$ contains a stack of length three,
implying $\vartheta \notin G_{i-1}(\lambda^{i-1}_{j-})$, which is impossible.
Therefore, we have the bijection
\begin{equation} \label{E:G-formula}
\beta\colon G_{i-1}(\lambda^{i-1}_{j-})\longrightarrow
T^c_{i-3}(\lambda^{i-3}_{j-})\setminus G_{i-3}(\lambda^{i-3}_{j-}),
\end{equation}
from which we conclude
\begin{equation*}
g_{i-1}(\lambda^{i-1}_{j-})=t^c_{i-3}(\lambda^{i-3}_{j-})-g_{i-3}
(\lambda^{i-3}_{j-}).
\end{equation*}
Replacing the term $g_r(\lambda^{r}_{j-})$ and using the fact that for
any shape $\mu$, $g_1(\mu)=g_0(\mu)=0$ holds, we arrive at
\begin{equation*}
g_{i-1}(\lambda^{i-1}_{j-})=
 \sum_{p=1}^{\lfloor \frac{i-1}{2}
\rfloor}(-1)^{p-1}t^c_{i-2p-1}(\lambda^{i-2p-1}_{j-}).
\end{equation*}
This allows us to rewrite eq.~(\ref{E:union}) as
\begin{equation*}
 t^c_i(\lambda^i)=
t^c_{i-1}(\lambda_0^{i-1})+
\sum_{j=1}^{k-1}t^c_{i-1}(\lambda^{i-1}_{j+})+
\sum_{j=1}^{k-1} \sum_{p=0}^{\lfloor \frac{i-1}{2}
\rfloor}(-1)^p  t^c_{i-1-2p}(\lambda^{i-1-2p}_{j-})
\end{equation*}
and the proof of the lemma is complete.
\qed

Lemma~\ref{L:recursion} allows us to compute the terms $t^c_i(\lambda)$ for
arbitrary $i$ and $\lambda$ recursively via the terms $t^c_h(\lambda')$, where
$h<i$ and the shapes $\lambda'$ differ from $\lambda$ by at most
one square.

We next generate a $*$-tableaux $\vartheta \in T^c_n(\lambda^n=\varnothing)$ from
right to left. For this purpose we set $\mu^i=\lambda^{n-i}$ for all $0\le i \le n$
and initialize $\mu^0=\varnothing$. Suppose we have at step $i$ the shape $\mu^i$ and
consider the $T^c_{n-i}(\lambda^{n-i})$-paths starting from $\lambda^0=\varnothing$
and ending at $\lambda^{n-i}=\mu^i$.

\begin{corollary}\label{C:oe}
The transition probabilities
\begin{equation} \label{E:transp}
\mathbb{P}(X^{i+1}=\mu^{i+1}\mid X^{i}=\mu^{i})=
\begin{cases}
 \frac{t^c_{n-i-1}(\mu^{i+1})}{t^c_{n-i}(\mu^{i})} &
\mu^{i}\setminus\mu^{i+1}=+\square_j, \varnothing \\
 \frac{\sum_{p=0}^{\lfloor (n-i-1)/2 \rfloor} (-1)^p
t^c_{n-i-2p-1}(\mu^{i+1})}{t^c_{n-i}(\mu^{i})} &
\mu^{i}\setminus \mu^{i+1}=-\square_j,
\end{cases}
\end{equation}
where $1\le j\le k-1$, induce a locally uniform Markov-process $(X^i)_i$ whose
sampling paths are shape-sequences $(\mu^i)_i$.
\end{corollary}

Let {\bf Rand}$(\mu^{i})$ denote the random process of locally
uniformly choosing $X^{i+1}=\mu^{i+1}$ for given $X^{i}=\mu^{i}$
using the transition probabilities given in eq.~(\ref{E:transp}).
Corollary~\ref{C:oe} gives rise to the following algorithm:

\begin{algorithm}
\caption{{\bf Core}$(n,k)$} \label{A:core}
\begin{algorithmic} [1]
\STATE $m \leftarrow 0$
\WHILE {$m<n$}
\STATE $\mu^{m+1} \leftarrow \textbf{Rand}(\mu^{m})$
\IF {$\mu^{m+1}\setminus \mu^{m}=+\square$}
\STATE insert $(m+1)$ in the new square
\ELSIF {$\mu^{m+1}\setminus \mu^{m}=-\square$}
\STATE let $pop$ be the unique extracted entry of $T^{m}$
which if RSK-inserted into $T^{m+1}$ recovers $T^{m}$
\STATE create an arc $(pop,m+1)$
\IF {$(pop,m+1)$ is stacking with $lastpair$}
\STATE restart the process {\bf Core}$(n,k)$
\ELSE
\STATE put $(pop,m+1)$ in the arc set $A$
\STATE $lastpair \leftarrow (pop,m+1)$
\ENDIF
\ENDIF
\STATE $m \leftarrow m+1$
\ENDWHILE
\end{algorithmic}
\end{algorithm}

The key observation now is that any core-diagram generated via the
above Markov process has {\it uniform} probability.

\begin{theorem}\label{T:eins}
Any core-diagram generated via the Markov-process $(X^i)_i$ (by means of
the algorithm ${\bf Rand}(\mu^{i})$) is generated with uniform probability.
\end{theorem}

\proof Suppose we are given a sequence of shapes
$$
\mu^i,\mu^{i-1},\ldots,\mu^0=\varnothing
$$
Let $U_{n-i}(\mu^{i})$ denote the subset of $*$-tableaux
\begin{equation*}
\varnothing=\lambda^0,\lambda^{1},\ldots, \lambda^{n-i}=\mu^i
\end{equation*}
such that there is no stack in the induced arc set of
\begin{equation*}
(\lambda^0,\ldots,
\lambda^{n-i-1},\lambda^{n-i}=\mu^i,\mu^{i-1},\ldots,\mu^0=\varnothing).
\end{equation*}
In particular, $U_n(\varnothing)$ denotes the set of all $*$-tableaux of
shape $\varnothing$ having at most $(k-1)$ rows that generate only
core-diagrams.
Let $u_n(\varnothing)=|U_n(\varnothing)|$ denote the number of cores of
length $n$. By construction, we have
\begin{equation*}
U_{n-i}(\mu^{i}) \subseteq T^c_{n-i}(\mu^{i}),
\end{equation*}
We now condition the process $(X^{i})_i$, whose transition probabilities
are given by eq.~(\ref{E:transp}), on generating cores. That is, we consider
only those $*$-tableaux generated by $(X^{i})_i$ that are contained in
$U_n(\varnothing)$. Let this process be denoted by $(Z^{i})_i$. We
observe
\begin{eqnarray*}
(T^c_{n-i-1}(\mu^{i+1}) \setminus G_{n-i-1}(\mu^{i+1}))
\, \cap \, U_{n-i-1}(\mu^{i+1}) & = & U_{n-i-1}(\mu^{i+1}) \\
T^c_{n-i}(\mu^{i})\, \cap \, U_{n-i}(\mu^{i}) & = &
U_{n-i}(\mu^{i}) \\
T^c_{n-i-1}(\mu^{i+1})\, \cap\, U_{n-i-1}(\mu^{i+1}) & = &
U_{n-i-1}(\mu^{i+1}).
\end{eqnarray*}
Accordingly, using eq.~(\ref{E:transp}), we derive for the transition
probabilities
\begin{eqnarray*}
\mathbb{P}(Z^{i+1}\mid Z^{i}) &=&
\frac{|U_{n-i-1}(\mu^{i+1})|}{|U_{n-i}(\mu^{i})|}.
\end{eqnarray*}
Therefore we arrive at
\begin{equation*}
\mathbb{P}(Z^{i+1}
) = \prod_{p=0}^{i}
      \frac{|U_{n-i-1+p}(\mu^{i+1-p})|}{|U_{n-i+p}(\mu^{i-p})|}
       = \frac{|U_{n-i-1}(\mu^{i+1})|}{|U_n(\mu^0=\varnothing)|}
       = \frac{|U_{n-i-1}(\mu^{i+1})|}{u_n(\varnothing)}
\end{equation*}
and in particular
\begin{equation*}
\mathbb{P}(Z^n=\varnothing)=\frac{|U_0(\mu^n=\varnothing)|}
{|U_n(\mu^0=\varnothing)|}= \frac{1}{u_n(\varnothing)},
\end{equation*}
which implies that the process $(Z^i)_i$ generates cores with uniform
probability. \qed

\section{Modular diagrams}\label{S:can}

Any $\sigma$-modular diagram can be mapped into a $\sigma$-weighted core, i.e.~a
diagram whose arcs have additional weights $\ge \sigma$. Suppose we
have a $*$-tableaux of $\varnothing$, $\vartheta$, whose induced diagram is a
$\sigma$-modular diagram. Repeated application of Lemma~\ref{L:removal} for 
each respective stack 
$$
S=\left((p,q),(p+1,q-1),\ldots,(p+(s-1),q-(s-1))\right),
$$
allows us to replace any insertion-step $p+1,\ldots, p+(s-1)$ as well as any 
extraction-step $q-(s-1),\ldots, q-1$ by $\varnothing$-steps, respectively.
Removing the $2(s-1)$ $\varnothing$-steps and assigning the stack-lengths $s$
to the {\it extraction} in step $q$, generates a $*$-tableaux of $\varnothing$
with weights, $\theta$ ($\sigma$-weighted $*$-tableaux).

\begin{figure}[ht]
\centerline{\epsfig{file=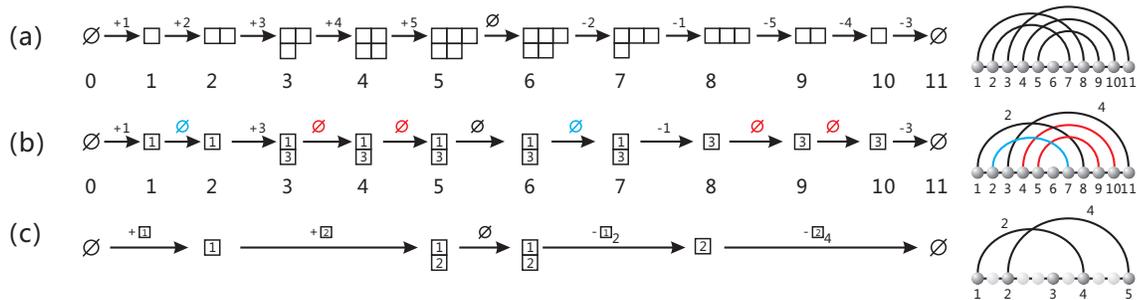,width=1\textwidth} \hskip8pt}
\caption{\small (a) a $*$-tableaux whose induced diagram is a $2$-modular
diagram. 
(b) the $*$-tableaux obtained by repeated application of Lemma~\ref{L:removal}. 
The red and blue removed arcs correspond the red and blue $\varnothing$-steps 
in the $*$-tableaux, respectively. 
(c) the weighted $*$-tableaux induced by (b) with weights $2$ and $4$ assigned 
to the two extraction steps, respectively, and its induced weighted core.
} \label{F:seq_can}
\end{figure}

Using the correspondence between $*$-tableaux and diagrams, a $\sigma$-weighted
core can therefore be represented as a sequence of shapes, $\theta$ in which,
preceding each extraction step, we have the additional insertion of exactly
$2(s-1)$ $\varnothing$-steps, see Fig.~\ref{F:seq_can}.
Let $W^\sigma_i(\lambda^r)$ denote the set of $\sigma$-weighted $*$-tableaux.
Each such $\theta \in W^\sigma_i(\lambda^r)$ induces a unique $*$-tableaux,
$p(\theta)$, contained in $T^c_r(\lambda^r)$ and we have
\begin{equation*}
i=r+\sum_{\ell=1}^{h\le r/2}2(s_\ell-1),
\end{equation*}
where $s_\ell$ is the weight of the $\ell$th extraction in $\theta$.
We set $w^\sigma_i(\lambda^r)=|W^\sigma_i(\lambda^r)|.$

\begin{lemma}\label{L:rec-w}
We have the recursion formula
\begin{eqnarray}
w_i(\lambda^r) & = & w^\sigma_{i-1}(\lambda_0^{r-1})+
                     \sum_{j=1}^{k-1}w_{i-1}^\sigma(\lambda_{j+}^{r-1}) \\
& & +\sum_{j=1}^{k-1}\sum_{s=\sigma}^{\lfloor \frac{i+1}{2}\rfloor}
\sum_{\ell=1}^{\lfloor \frac{s}{\sigma}
\rfloor}(-1)^{\ell-1} p(s,\ell,\sigma)w^\sigma_{i-2s+1}(\lambda_{j-}^{r-1}),
\nonumber
\end{eqnarray}
where $p(a,\ell,\sigma)$ denotes the number of partitions of $a$ into
$\ell$ blocks, $\{a_1,a_2,\ldots,a_{\ell}\}$, such that $\forall i\le \ell$,
$a_i\ge \sigma$.
\end{lemma}
\proof
Any $*$-tableaux $\theta \in W_i^\sigma(\lambda^r)$, where $i=r+\sum_{\ell=1}^h
2(s_\ell-1)$, $s_\ell$ is the weight assigned to the $\ell$th extraction step in
$\theta$.
We consider the weighted $*$-tableaux, $\theta^\prime$, derived from $\theta$ by
removing the shape in step $r$. If $\lambda^r$ is derived from $\lambda^{r-1}$
by doing nothing, then $\theta^\prime \in W^\sigma_{i-1}(\lambda^{r-1}_0)$.
Similarly, if $\lambda^r$ is derived from $\lambda^{r-1}$ by adding a square in
the $j$th row, we have $\theta^\prime \in W^\sigma_{i-1}(\lambda^{r-1}_{j+})$.
In case of $\lambda^r$ being derived from $\lambda^{r-1}$ via removing a square
from the $j$th row, we are given an extraction step with associated weight $s$.
Thus,
\begin{equation*}
\theta^\prime \in
W^\sigma_{r-1+\sum_{\ell=1}^{h-1}2(s_\ell-1)}(\lambda^{r-1}_{j-})
=W^\sigma_{(i-1)-2(s-1)}(\lambda^{r-1}_{j-}).
\end{equation*}
$\theta^\prime$ determines $*$-tableaux, $p(\theta')\in
T^c_{r-1}(\lambda^{r-1}_{j-})
\setminus G_{r-1}(\lambda^{r-1}_{j-})$.
Let $V^\sigma_{(i-1)-2(s-1)}(\lambda^{r-1}_{j-})$ denote the set
of weighted $*$-tableaux $\theta_1$ such that
$p(\theta_1)\in G_{r-1}(\lambda^{r-1}_{j-})$. We set $v^\sigma_{(i-1)-2(s-1)}
(\lambda^{r-1}_{j-})=|V^\sigma_{(i-1)-2(s-1)}(\lambda^{r-1}_{j-})|$.
Note that then $\theta^\prime \in
W^\sigma_{(i-1)-2(s-1)} (\lambda^{r-1}_{j-}) \setminus
V^\sigma_{(i-1)-2(s-1)}(\lambda^{r-1}_{j-})$, whence
\begin{equation*}
\begin{split}
W_i^\sigma(\lambda^r)&=W_{i-1}^\sigma(\lambda_0^{r-1}) \,\dot\cup
 \left(\bigcup_{j=1}^{k-1} W_{i-1}^\sigma(\lambda^{r-1}_{j+})\right) \\
 &\dot\cup
\left(\bigcup_{j=1}^{k-1}\bigcup_{s=\sigma}^{\lfloor \frac{i+1}{2} \rfloor}
[W_{(i-1)-2(s-1)}^\sigma(\lambda^{r-1}_{j-}
)\setminus V^\sigma_{(i-1)-2(s-1)}(\lambda^{r-1}_{j-})]\right).
\end{split}
\end{equation*}
We therefore derive
\begin{equation} \label{E:W-formula}
w_i^\sigma(\lambda^r)= w^\sigma_{i-1}(\lambda^{r-1}_0)+
             \sum_{j=1}^{k-1}w^\sigma_{i-1}(\lambda^{r-1}_{j+})+
             \sum_{j=1}^{k-1}\sum_{s=\sigma}^{\lfloor \frac{i+1}{2} \rfloor}
             [w^\sigma_{i-2s+1}(\lambda^{r-1}_{j-})-v^\sigma_{i-2s+1}
             (\lambda^{r-1}_{j-})].
\end{equation}
We proceed by considering a $*$-tableaux $\zeta \in V^\sigma_{(i-1)-2(s-1)}
(\lambda^{r-1}_{j-})$. By construction, $(r-1)$ is a $-\square$-step.
Suppose the induced arc of this extraction is $\alpha$ and the weight
assigned to it is given by $s^\prime$. Then $p(\zeta)\in G_{r-1}
(\lambda^{r-1}_{j-})$ and we have the bijection
\begin{equation*}
\beta\colon G_{r-1}(\lambda^{r-1}_{j-})\longrightarrow
T^c_{r-3}(\lambda^{r-3}_{j-})\setminus G_{r-3}(\lambda^{r-3}_{j-}),
\end{equation*}
obtained by removing the insertion and extraction step of the extracted square
in step $(r-1)$. Taking into the account weights, $\beta$ gives rise
to the bijection
\begin{equation*}
\beta^\prime \colon V^\sigma_{(i-1)-2(s-1)}(\lambda^{r-1}_{j-})\longrightarrow
\bigcup_{s^\prime=\sigma}^{\lfloor \frac{i-2s+1}{2} \rfloor}
[W^\sigma_{(i-1)-2(s+s^\prime-1)}(\lambda^{r-3}_{j-})\setminus
V^\sigma_{(i-1)-2(s+s^\prime-1)}(\lambda^{r-3}_{j-})],
\end{equation*}
from which we conclude
\begin{equation*}
 v^\sigma_{(i-1)-2(s-1)}(\lambda^{r-1}_{j-}) = \sum_{s^\prime=\sigma}^{\lfloor
\frac{i-2s+1}{2}
\rfloor}[w_{(i-1)-2(s+s^\prime-1)}^\sigma(\lambda^{r-3}_{j-}
)-v^\sigma_{(i-1)-2(s+s^\prime-1)}(\lambda^{r-3}_{j-})].
\end{equation*}
Using (a)
$$
\sum_{s_1=\sigma}\ldots \sum_{s_\ell=\sigma}x_{s_1+\cdots +s_\ell}=
p(s,\ell,\sigma)x_s
$$
where $p(s,\ell,\sigma)$ denotes the number of partitions of $s$ into
$\ell$ blocks of size $\ge \sigma$, and (b) that for any shape $\mu$,
$v^\sigma_1(\mu)=v^\sigma_0(\mu)=0$. We iterate the above formula by
replacing the terms $v^\sigma_r(\lambda^r_{j-})$
\begin{equation*}
\begin{split}
& \sum_{s=\sigma}^{\lfloor \frac{i+1}{2}\rfloor}(w^\sigma_{i-2s+1}
(\lambda^{r-1}_{j-})-v^\sigma_{i-2s+1}(\lambda^{r-1}_{j-})) \\
=& \sum_{s=\sigma}^{\lfloor \frac{i+1}{2}\rfloor}w^\sigma_{i-2s+1}
(\lambda^{r-1}_{j-})
-\sum_{s=\sigma}^{\lfloor \frac{i+1}{2}\rfloor}\sum_{s^\prime=\sigma}^{\lfloor
\frac{i-2s+1}{2}\rfloor}(w^\sigma_{i-2s-2s^\prime+1}
(\lambda^{r-3}_{j-})-v^\sigma_{i-2s-2s^\prime+1}(\lambda^{r-3}
_{j-})) \\
\vdots \\
=&\sum_{s_1=\sigma}\ldots\sum_{s_\ell=\sigma}^{2(s_1+\cdots+s_\ell)\le
i+1} (-1)^\ell w^\sigma_{i-2(s_1+\cdots+s_\ell)+1
}(\lambda^{r-2\ell+1}_{j-}) \\
=& \sum_{s=\sigma}^{\lfloor \frac{i+1}{2}\rfloor} \sum_{\ell=1}^{\lfloor
\frac{s}{\sigma}
\rfloor}(-1)^{\ell-1} p(s,\ell,\sigma)w^\sigma_{i-2s+1}(\lambda^{r-1}_{j-}),
\end{split}
\end{equation*}
whence the lemma.  \qed

Lemma~\ref{L:rec-w} allows us to compute $w^\sigma_{i}(\mu)$ for arbitrary
$i,\mu$ inductively via the terms $w^\sigma_h(\lambda)$ and $h<i$.
We next consider the generation of a $*$-tableaux, $\vartheta$, which corresponds
to a $\sigma$-modular diagram.
For this purpose we shall generate a weighted $*$-tableaux
$\theta\in W^\sigma_n(\lambda^m=\varnothing)$.
Taking the sum over all weights we have
$m=n-\sum_h 2(s_h-1)$.
We construct $\theta$ inductively from right to left setting
$\mu^r=\lambda^{m-r}$.
We initialize $\mu^0=\varnothing$ and assign in case of $\mu^r\setminus
\mu^{r+1}=-\square$ a weight to step $r$. Suppose we have arrived at
$\mu^{r}$, with the corresponding set of weights, $S^r$. Considering sequences
of weighted $*$-tableaux contained in $W^\sigma_{n-r-\sum_{s_\ell \in
S^{r}}2(s_\ell-1)}(\mu^{r})$, Lemma~\ref{L:rec-w} implies

\begin{corollary}\label{C:kkk}
Let
\begin{eqnarray*}
t & = & n-r-\sum_{s_\ell\in S^r}2(s_\ell-1) \\
z^\sigma_{t-1}(\mu^{r+1}_{j-},s) & = & \sum_{\ell=1}^{\lfloor
\frac{s}{\sigma}\rfloor}(-1)^{\ell-1}
p(s,\ell,\sigma)w^\sigma_{t-2s+1}(\mu_{j-}^{r+1}).
\end{eqnarray*}
The transition probabilities
\begin{equation} \label{E:trans_can}
\begin{split}
& \mathbb{P}(X^{r+1}=(\mu^{r+1},S^{r+1})\mid X^r=(\mu^{r},S^{r})) \\
=& \begin{cases}
 \frac{w^\sigma_{t-1}(\mu^{r+1})}{w^\sigma_{t}(\mu^{r})} &
 \mu^{r}\setminus \mu^{r+1}= +\square_j, S^{r+1}=S^{r} \\
 \frac{w^\sigma_{t-1}(\mu^{r+1})}{w^\sigma_{t}(\mu^{r})} &
\mu^{r}\setminus\mu^{r+1}= \varnothing, S^{r+1}=S^{r} \\
 \frac{z^\sigma_{t-1}(\mu^{r+1},s)}{w^\sigma_{t}(\mu^{r})} &
\mu^{r}\setminus\mu^{r+1}= -\square_j, S^{r+1}=S^{r}\cup \{s\},
\end{cases}
\end{split}
\end{equation}
for $1\le j \le k-1$, generate a locally uniform Markov-process $(X^i)_i$.
\end{corollary}

Corollary~\ref{C:kkk} represents an algorithm for constructing $\sigma$-modular
diagrams. In analogy to the case of core-diagrams, if $X$ successfully
constructs a modular diagram, it generates the latter with uniform probability,
see Fig.~\ref{F:dis_suc} (lefthand side).
\begin{algorithm}
\caption{{\bf Canonical}$(n,k,\sigma)$} \label{A:canonical}
\begin{algorithmic} [1]
\STATE $m \leftarrow 0$
\WHILE {$m<n$}
\STATE $(\mu^{m+1},size) \leftarrow \textbf{RandStep}(\mu^{m})$
\IF    {$\mu^{m+1} \setminus \mu^{m}=+\square$}
\STATE insert $(m+1)$ in the new square
\STATE assign $size$ to the the new square
\STATE $m\leftarrow m+size-1$
\ELSIF {$\mu^{m+1} \setminus \mu^{m} =-\square$}
\STATE let $pop$ be the unique extracted entry of $T^{m}$
which if RSK-inserted into $T^{m+1}$, recovers $T^m$ and let
$size$ be the integer assigned to the extracted square 
\STATE create a stack $\{(pop,m+size),\cdots,(pop+size-1,m+1)\}$
\IF {$(pop,m+size)$ is stacking with $lastpair$}
\STATE restart the process {\bf Canonical}$(n,k,\sigma)$
\ELSE
\STATE put $\{(pop,m+size),\cdots,(pop+size-1,m+1)\}$ in the arc set $A$
\STATE $lastpair \leftarrow (pop,m+size)$
\STATE $m\leftarrow m+size-1$
\ENDIF
\ENDIF
\STATE $m \leftarrow m+1$
\ENDWHILE
\end{algorithmic}
\end{algorithm}
Consequently, the process $(X^i)_n$ generates random $\sigma$-modular,
$k$-noncrossing diagram in $O(n)$ time and space complexity. According to
the recursion of Lemma~\ref{L:rec-w}, we compute $w^\sigma_i(\lambda^i)$ for
arbitrary $\lambda^i$ with at most $(k-1)$ rows and all $i\le n$ in
$O(n)\times O(n^{k-1})=O(n^k)$ time and space complexity.
\begin{theorem}\label{T:zwei}
Any modular diagram derived via the Markov-process $(X^r)_r$ is
generated with uniform probability.
\end{theorem}

\proof Suppose we have a sequence of shapes
$$
\mu^{r},\mu^{r-1},\ldots, \mu^0=\varnothing
$$
with weights assigned to each $\mu^{i-1}\setminus \mu^i=-\square_j$-step
and set of weights, $S^r$.
Let
$$
t = n-r-\sum_{s_\ell\in S^{r}}2(s_\ell-1)
$$
and $D^\sigma_{m-r}(\mu^{r})$ be the set of weighted $*$-tableaux
\begin{equation*}
\lambda^0=\varnothing,\lambda^{1},\ldots, \lambda^{m-r}=\mu^r
\end{equation*}
such that
\begin{equation*}
\varnothing=\lambda^0,\ldots,\lambda^{m-r-1},\lambda^{m-r}=\mu^r,\mu^{r-1},
\ldots,\mu^0=
\varnothing
\end{equation*}
is contained in $W_n^\sigma(\lambda^m=\varnothing)$.
Summing over all weights we have $m=n-\sum_{s_h}2(s_h-1)$ and by construction,
$D^\sigma_{m-r}(\mu^r)\subseteq W^\sigma_{t}(\mu^r)$.
In particular
\begin{equation*}
d^\sigma_m(\varnothing)=|D^\sigma_m(\mu^m=\varnothing)|
\end{equation*}
equals the number of weighted core of length $m$, i.e.~the number of
modular diagrams of length $n$. Suppose now we only consider sampling
paths of weighed cores generated via $(X^r)_r$ (whose transition probabilities
is given by eq.~(\ref{E:trans_can})) contained in $D^\sigma_m(\varnothing)$.
We denote the resulting process by $(Z^r)_r$. In view of
\begin{equation*}
D^\sigma_{m-r}(\mu^r)\subseteq W^\sigma_{t}(\mu^r),
\end{equation*}
we observe that
\begin{eqnarray*}
(W^\sigma_{t-1}(\mu^{r+1}) \setminus V^\sigma_{t-1}(\mu^{r+1})) \cap
D^\sigma_{m-r-1}(\mu^{r+1}) & = & D^\sigma_{m-r-1}(\mu^{r+1}) \\
W^\sigma_{t}(\mu^{r}) \cap D^\sigma_{m-r}(\mu^{r}) & = &
D^\sigma_{m-r}(\mu^{r}) \\
W^\sigma_{t-1}(\mu^{r+1}) \cap D^\sigma_{m-r-1}(\mu^{r+1}) & = &
D^\sigma_{m-r-1}(\mu^{r+1}).
\end{eqnarray*}
Therefore we have
$$
\mathbb{P}(Z^{r+1}\mid Z^r)=\frac{|D_{m-r-1}^\sigma(\mu^{r+1})|}
{|D_{m-r}^\sigma(\mu^r)|}
$$
and consequently
\begin{equation}
\mathbb{P}(Z^{r+1})=\prod_{p=0}^{r}\frac{|D_{m-r-1+p}^\sigma(\mu^{r+1-p})|}
{|D_{m-r+p}^\sigma(\mu^{r-p})|}=\frac{|D_{m-r-1}^\sigma(\mu^{r+1})|}
{|D_{m}^\sigma(\mu^{0}=\varnothing)|}.
\end{equation}
In particular,
$$
\mathbb{P}(Z^m=\varnothing)=\frac{|D^\sigma_0(\mu^m=\varnothing)|}
{|D_{m}^\sigma(\mu^{0}=\varnothing)|}=\frac{1}{|D_{m}^\sigma(\mu^{0}
=\varnothing)|}=\frac{1}{d^\sigma_m(\varnothing)}.
$$
That is, the process $(Z^r)$ generates modular diagrams uniformly and the
theorem follows.\qed

\begin{figure}[ht]
\centerline{\epsfig{file=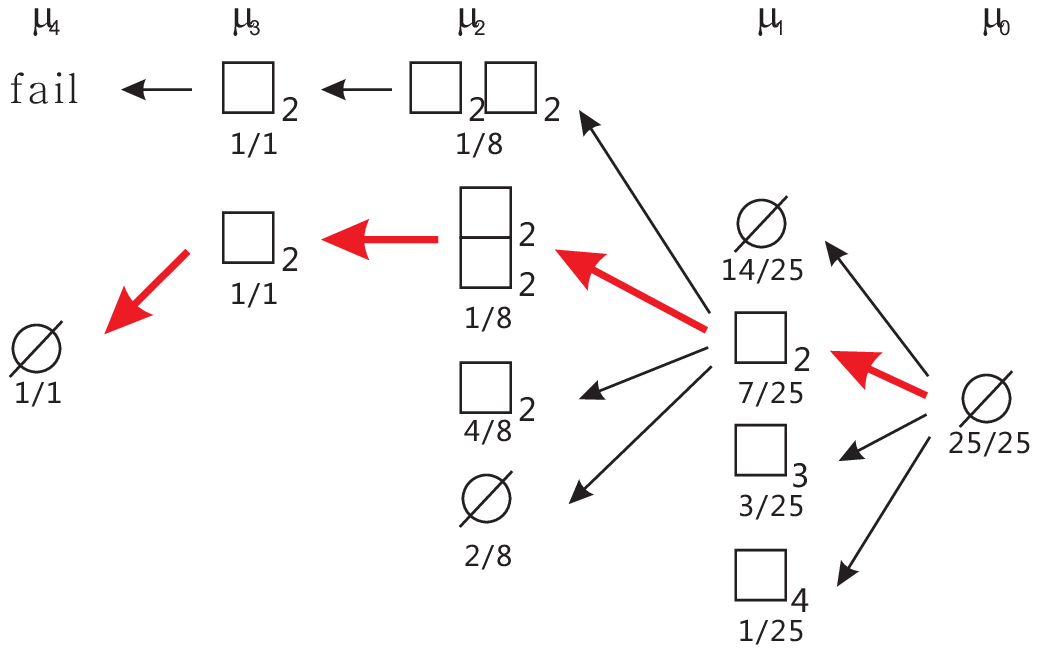,width=0.5\textwidth} \hskip8pt}
\caption{\small Building weighted $*$-tableaux via the transition probabilities
given in eq.~(\ref{E:trans_can}).
Here, the top path fails to generate a $2$-modular diagram while the red path
succeeds. According to Theorem~\ref{T:zwei} each such modular diagram is
generated with uniform probability.} \label{F:trans_seq}
\end{figure}

\begin{figure}[ht]
\centerline{\epsfig{file=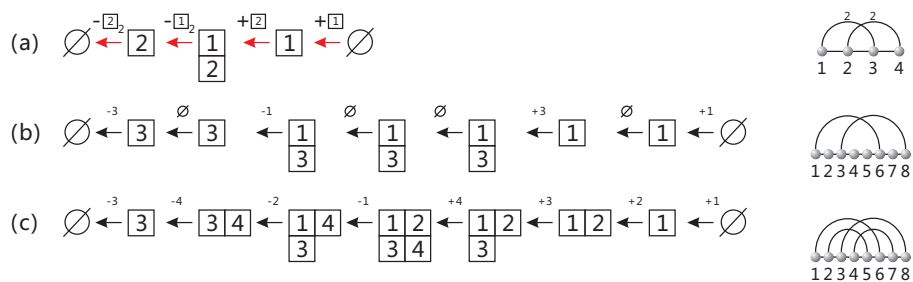,width=0.8\textwidth} \hskip8pt}
\caption{\small (a) the red path of Fig.~\ref{F:trans_seq}.
(b) the $*$-tableaux derived by adding four
$\varnothing$-steps in (a). (c) adding two pairs of insertion and extraction
steps, which produces the $2$-modular diagram.
} \label{F:rebuild}
\end{figure}

In Fig.~\ref{F:trans_seq}, we showcase two paths constructed via the
process {\bf Canonical}$(n,k,\sigma)$, for $n=8, k=3$ and $\sigma=2$.
In Fig.~\ref{F:rebuild} we construct the corresponding $2$-modular
diagram from the red path displayed in Fig.~\ref{F:trans_seq}.

\bibliographystyle{plain}  

\end{document}